\documentclass[12pt]{amsart}
\usepackage{graphicx,amssymb,latexsym,bm,color}
\usepackage[dvips,centering,includehead,width=15.1cm,height=22.7cm]{geometry}
\usepackage{amsmath}
\usepackage{graphicx}
\usepackage{epsfig}
\usepackage{amssymb}
\usepackage{subfigure}
\input{xy} \xyoption{all}

\newtheorem{theorem}{Theorem}[section]

\theoremstyle{definition} \newtheorem{definition}[theorem]{Definition}
\newtheorem{example}[theorem]{Example}
\newtheorem{remark}[theorem]{Remark}

\newcommand{\eps}{{\varepsilon}}
 \newcommand{\N}{{\mathcal N}}
\newcommand{\C}{{\mathcal C}} 
\renewcommand{\L}{{\mathcal L}}

 \renewcommand{\P}{{\mathcal P}}
\renewcommand{\O}{{\mathcal O}}

\DeclareMathOperator{\wt}{wt}

\DeclareMathOperator{\trop}{trop}
\DeclareMathOperator{\primitive}{primitive}
\DeclareMathOperator{\mult}{mult}
\DeclareMathOperator{\Log}{Log}
\DeclareMathOperator{\Aut}{Aut}
\DeclareMathOperator{\conv}{conv}

\DeclareMathOperator{\id}{id}
\DeclareMathOperator{\abs}{abs}
\DeclareMathOperator{\interior}{int}

\newcommand{\RR}{{\mathbb R}}
\newcommand{\ZZ}{{\mathbb Z}}

\newcommand{\QQ}{{\mathbb Q}}
\newcommand{\PP}{{\mathbb C}{\mathbb P}}
\newcommand{\CC}{{\mathbb C}}

\title[Counting Algebraic Curves with Tropical Geometry]{Counting Algebraic Curves \\ with Tropical Geometry} 
\author{Florian Block}
\date{\today}
\address{Florian Block, Mathematics Institute, University of Warwick,
Coventry, CV4 7AL,
United Kingdom.}
\email{f.s.block@warwick.ac.uk}
\thanks{The author was supported by
the EPSRC grant EP/I008071/1.}

\begin{document}

\begin{abstract}
Tropical geometry is a piecewise linear ``shadow'' of algebraic
geometry. It allows for the computation of several cohomological
invariants of an algebraic variety. In particular, its application to
enumerative algebraic geometry led to significant progress.

In this survey, we give an introduction to tropical geometry techniques
for  algebraic curve counting problems. We also survey some recent developments,
with emphasis on the computation of the degree of
Severi varieties of the complex projective plane and other toric
surfaces as well as Hurwitz numbers and applications to real
enumerative geometry. This paper is based on the author's lecture
at the Workshop on Tropical Geometry and Integrable Systems in
Glasgow, July 2011.
\end{abstract}

\maketitle

\vspace{-0.8em}

\section{Enumerative Algebraic Geometry}
\label{sec:enum_AG}

\subsection{Overview}

Enumerative algebraic geometry is the study of enumerations of algebro-geometric
objects with certain properties. In this article, we mostly consider the
enumeration of complex algebraic curves. A
typical question is: ``What is the number $N_{d, 0}$ of irreducible rational curves in the complex
plane $\PP^2$ of degree
$d$ passing through $3d-1$ points in general position?'' 

Enumeration of algebraic curves in a given algebraic variety $X$ is closely
related to its Gromov-Witten theory. If $X$ is a del Pezzo surface
(i.e., a projective algebraic surface with ample anticanonical bundle) its
Gromov-Witten invariants are enumerative, which means that they can be
computed by a curve enumeration~\cite{Va00}. 
For example,  the numbers $N_{d, 0}$ are
the \emph{rational Gromov-Witten invariants} of $X =  \PP^2$.
Classically, we have $N_{1,0} = N_{2, 0} = 1$ and $N_{3,0} = 12$.
In the late 19th century, Zeuthen computed $N_{4,0} = 620$. The
number $N_{5,0} = 87304$ was computed in the mid 20th century. 
For larger $d$, no progress was made until Kontsevich~\cite{KM94}, in
1995, computed $N_{d, 0}$, for all $d$, by his famous recursion
\[
N_{d,0} = \sum_{\substack{
d_1 + d_2 = d\\
d_1, d_2 >0
}} 
\left( d_1^2 d_2^2
  \binom{3d-4}{3d_1 - 2} - d_1^3 d_2 \binom{3d-4}{3d_1 - 1} \right)
N_{d_1, 0} N_{d_2, 0}.
\]

Using
tropical geometry, Gathmann and Markwig~\cite{GM053} reproved
Kontsevich's formula, based on Mikhalkin's Correspondence Theorem between
algebraic and tropical plane curves (see
Theorem~\ref{thm:correspondencetheorem}). An outline of their proof is
given in \cite[Section~3.2]{Ga06}.

More generally, one can allow curves of arbitrary genus and ask for the number
$N_{d, g}$ of irreducible degree-$d$ genus-$g$ plane curves passing
through $3d + g - 1$ points in general position. The numbers $N_{d, g}$ are the
\emph{Gromov-Witten invariants} of $ \PP^2$. The Gromov-Witten
invariants $N_{d, g}$ were computed by Caporaso  
and Harris~\cite{CH98} for all $d$ and $g$ in 1998. Even more generally, one can
consider appropriate counts of genus-$g$ curves on other surfaces, or
in algebraic varieties of higher dimension. The enumerative meaning,
however, can in general be quite subtle.

Closely related to the Gromov-Witten invariant $N_{d, g}$ is the
\emph{Severi degree} $N^{d, \delta}$, counting plane curves of degree 
$d$ with exactly $\delta$ nodes as singularities (we call such curves
\emph{$\delta$-nodal}) passing through
$\tfrac{(d+3)d}{2} - \delta$ points in $\PP^2$ in general position. Equivalently,
$N^{d, \delta}$ is the degree of the \emph{Severi variety}
parametrizing such curves. Enrique~\cite{En12} and Severi~\cite{Se21}
introduced these varieties around 100 years ago.

 Tropical geometry techniques have been applied
successfully also to problems in \emph{real} enumerative geometry.
Later in this article, in Section~\ref{sec:Welschinger}, we briefly
mention \emph{Welschinger invariants}, a real analog of rational
Gromov-Witten invariants, and how they can be computed by tropical means.

Enumerative
algebraic geometry includes many further subjects, such as Schubert
Calculus. There, one considers questions of the form
``How many lines in $ \PP^3$
simultaneously intersect four given generic lines?'' (The answer is,
maybe surprisingly, two.) More generally, one counts linear subspaces,
or flags of subspaces, that meet given linear subspaces in a
prescribed way. One may expect that some of these question can also be
answered tropically in the future~\cite{SS04}.



\subsection{Enumerative Geometry on Toric Surfaces}

We now generalize the definitions of $N_{d, g}$ and $N^{d, \delta}$ to
toric surfaces. Such invariants can still be computed solely in terms
of tropical geometry (see
Theorem~\ref{thm:toriccorrespondence}). In Section~\ref{sec:toricpolys}, we
discuss an application of the resulting combinatorics implying
polynomiality of the curve counts, in some parameters of the
surface. This may suggest a generalization of the G\"ottsche
conjecture~\cite[Conjecture~2.1]{Go} to a family of possibly
non-smooth surfaces.

Fix a lattice polygon $\Delta$ in $\RR^2$, i.e., $\Delta$ is the convex hull of a
finite subset of $\ZZ^2$. As is well-known in toric
geometry, $\Delta$  determines, via its normal fan, a projective toric variety
$X = X(\Delta)$, together with an ample line bundle $\L = \L(\Delta)$ on
$X(\Delta)$. Conversely, any such data $(X, \L)$ determines a lattice
polygon. A common theme in toric geometry is that many geometric
invariants of $X(\Delta)$, such as smoothness or its Chow groups, can
be directly read off from the combinatorics of $\Delta$. For a detailed
introduction to toric varieties see Cox, Little, and Schenk's
recent book~\cite{CLS11} or Fulton's classical introduction~\cite{Fu93}.

Counting curves on $X(\Delta)$ of a given ``degree'' means counting
curves in the complete linear system $|\L(\Delta)|$ of $\L(\Delta)$
(or a subsystem thereof).
A concrete way to think about a curve $C$ in $|\L(\Delta)|$ is as
follows.
Let $f$ be a
polynomial (or Laurent polynomial) with Newton polygon $\Delta$. Then the closure, in
$X(\Delta)$, of the vanishing set of $f$ in the complex torus
$(\CC^*)^2$ is an element of $|\L(\Delta)|$.

Given a lattice polygon $\Delta$, the \emph{Severi degree}
$N^{\Delta, \delta}$ of the 
toric surface $X(\Delta)$, together with the line bundle $\L(\Delta)$, is
the number of (not necessarily irreducible) $\delta$-nodal curves in
$|\L(\Delta)|$ passing through $|\Delta \cap \ZZ^2| - 1 - \delta$
points in general position.

\begin{example}
Let $\Delta = \conv\{ (0,0),
(3,0), (0,2), (3,2)\}$ be the lattice polygon shown on the left of
Figure~\ref{fig:P1timesP1tropicalcurve}. Then $\Delta$ defines the toric surfaces $X(\Delta) =
\PP^1 \times \PP^1$ and the line bundle $\L(\Delta)$ equals $\O(3,2)$. The
elements of the linear system $|\O(3,2)|$ are the 
divisors in $\PP^1
\times \PP^1$ of 
polynomials in $x_0, x_1, y_0$, and $y_1$ of bi-degree $(3,2)$, where
$x_i$ and $y_i$ have degree $(1,0)$ and $(0,1)$, respectively. Thus,
$N^{\Delta, \delta}$ counts $\delta$-nodal curves in $\PP^1 \times
\PP^1$ of bi-degree $(3,2)$ through $(3+1)(2+1)-1-\delta = 11 - \delta$ points in
general position.
\end{example}

Notice that, unlike in the case of $\PP^2$, the number $N_{\Delta, g}$ of irreducible genus-$g$ curves
in $|\L(\Delta)|$ through sufficiently many points, in general, does
not equal a Gromov-Witten invariant of $X(\Delta)$.

\bigskip

{\bf Acknowledgement.} I thank the anonymous referee for valuable comments that
helped me to improve this paper.

\section{Tropical Geometry}
\label{sec:tropicalGeometry}

Tropical geometry is a piecewise linear analog (or ``shadow'') of algebraic
geometry. The main objects of study are \emph{tropical varieties}, i.e., weighted,
balanced, polyhedral complexes in a real vector space
$\RR^n$, equipped with a lattice $\ZZ^n \subset \RR^n$. We won't give the general definition here, and only discuss
the case of tropical curves in toric surfaces; see
Definition~\ref{def:tropicalcurve} for the $\PP^2$ case and
Definition~\ref{def:tropicalcurvedegreeDelta} for any toric
surface.

Introductory texts on tropical geometry include two book drafts, one by Maclagan
and Sturmfels~\cite{MS11}, the second
by Mikhalkin and Rau~\cite{MR11}. The former text is more extrinsic, with
tropical varieties often given as ``tropicalizations'' of 
algebraic varieties given by polynomial equations, and is more
computationally oriented.  The latter takes a more intrinsic
approach, with focus on developing  a theory of tropical
geometry in analogy with (non-tropical) algebraic geometry. For a
shorter introduction, with an emphasis on tropical curves, see
Gathmann's excellent survey~\cite{Ga06}.

\subsection{Tropical Curves for $\PP^2$}

\begin{definition}
\label{def:tropicalcurve}
A \emph{tropical plane curve of degree $d$} is a piecewise linear,
weighted graph $\Gamma$ in $\RR^2$ satisfying:
\begin{enumerate}
\item
\label{itm:weights}
all edges $e$ of $\Gamma$ have weights $\wt(e) \in \ZZ_{\ge 1}$,
\item
\label{itm:slopes}
all edges have rational slopes,
\item
\label{itm:directions}
the total weight of the edges of $\Gamma$ in each of the
  directions $\binom{-1}{0}$, $\binom{0}{-1}$, and $\binom{1}{1}$
  equals $d$, and $\Gamma$ has no other unbounded edges,
\item
\label{itm:balanced}
 all vertices $v$ of $\Gamma$ are \emph{balanced}, i.e.,
\vspace{-0mm}
\[
\sum_{\substack{
\text{edges e}\\
v \in e
}} 
\wt(e) \cdot \primitive(e, v) = 0,
\]
\vspace{-4mm}
\end{enumerate}
where $\primitive(e, v)$ is the \emph{primitive vector} of the edge
$e$ at the vertex $v$, i.e., the shortest non-zero integral vector in
the ray spanned by $e$.
\end{definition}

\begin{figure}[h]
\begin{center}
\includegraphics[height=4cm]{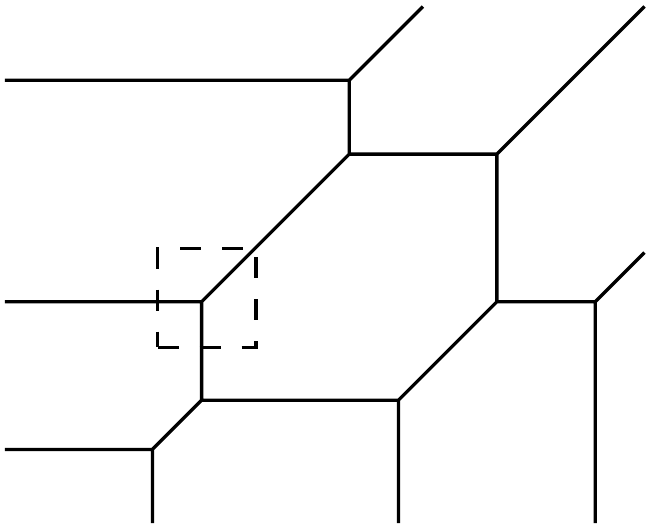}
\end{center}
\caption{A smooth tropical plane cubic. All edge weights are equal to
  $1$. A close-up of the highlighted  vertex is shown in
  Figure~\ref{fig:balancedvertex}.}
\label{fig:tropicalcubic}
\end{figure}

\begin{figure}[h]
\begin{center}
\includegraphics[height=3cm]{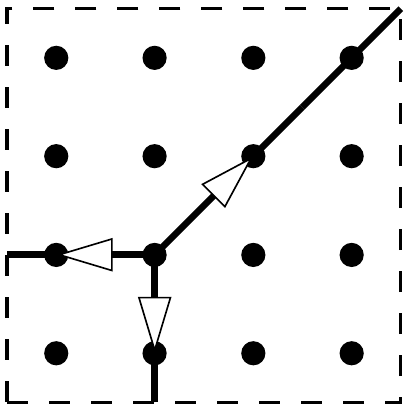}
\end{center}
\caption{A balanced vertex. The weighted sum of the adjacent primitive
  vectors vanishes:
  $1\binom{-1}{0} + 1\binom{0}{-1} + 1\binom{1}{1} =
  \binom{0}{0}$.}
\label{fig:balancedvertex}
\end{figure}

See Figure~\ref{fig:tropicalcubic} for an illustration of a  tropical plane curve of
degree $3$ and Figure~\ref{fig:balancedvertex} for a balanced vertex.
Condition~(\ref{itm:balanced}) in Definition~\ref{def:tropicalcurve}, also known as the \emph{zero-tension condition},
says that at each vertex a ``tug of war'', with directions given by the
edges and forces given by the weights, results in no net movement.

Every tropical plane curve $\Gamma$ has a number of numerical
invariants associated to it. We say that $\Gamma$ is
\emph{irreducible} if $\Gamma$ is not a union of two (non-empty)
tropical plane curves. If $\Gamma$ is irreducible, its \emph{genus}
$g(\Gamma)$ is the minimal first Betti number (or rank of the
fundamental group) of any topological graph $\Gamma'$ such that there
exists a surjective continuous map $\Gamma' \to \Gamma$. If $\Gamma$
has irreducible components $\Gamma_1, \dots, \Gamma_r$, then its genus
$g(\Gamma)$ is $\sum_{i=1}^r g(\Gamma_i) + 1 - r$.



The \emph{number of nodes} $\delta(\Gamma)$ of a tropical
plane curve of degree $d$ is $\tfrac{(d-1)(d-2)}{2} - g(\Gamma)$. This
formula is motivated by the corresponding genus-degree formula for
algebraic plane curves. Equivalently, if $\Gamma$ has irreducible
components of degree $d_1, d_2,
\dots$ and number of nodes  $\delta_1, \delta_2, \dots$, respectively, then
the number of nodes of $\Gamma$ equals $\sum_i 
\delta_i + \sum_{i < j} d_i d_j$. This formula parallels the
classical B\'ezout theorem for plane curves.
The interested reader wanting to ``find'' the locations of the
tropical nodes can have a look at~\cite[Theorem~2.9]{DT11}.

\begin{example}
The tropical plane degree-$3$ curve $\Gamma$ in Figure~\ref{fig:tropicalcubic}
has genus $g = 1$, the cycle is realized by the hexagon. The number of
nodes of $\Gamma$ is $\delta = \tfrac{(3-1)(2-1)}{2} - 1 = 0$. Thus
$\Gamma$ is a \emph{smooth tropical plane cubic}.
\end{example}

\begin{example}
Consider now the curve $\Gamma$ in Figure~\ref{fig:1nodalcubic}.
 It differs from the curve in Figure~\ref{fig:tropicalcubic} only around
the transverse crossing of two edges.
Although it looks like $\Gamma$ could have genus $1$, there is in fact
a surjective continuous map from a tree (considered as a topological
space) with nine leaves onto $\Gamma$. Thus $\Gamma$ has genus $0$,
and we have $\delta(\Gamma) = \tfrac{(3-1)(3-2)}{2} - 0 = 1$. (Here, the ``tropical
node'' is at the transverse intersection of the two edges.) Thus $\Gamma$ is
a \emph{rational tropical plane cubic}.


\begin{figure}[h]
\begin{center}
\includegraphics[height=4cm]{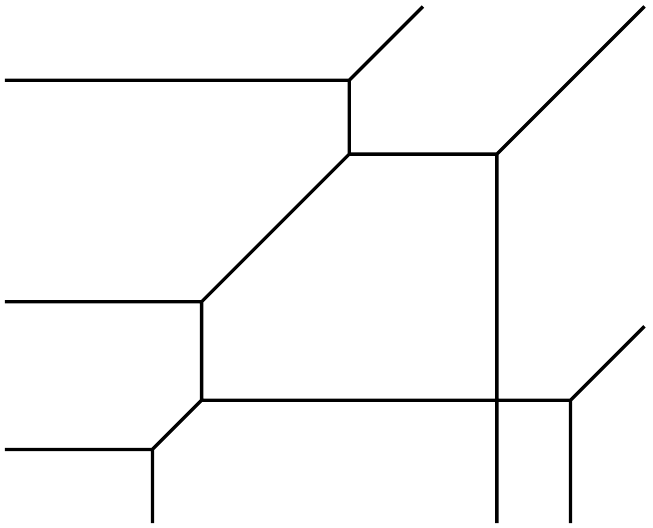}
\end{center}
\caption{A $1$-nodal tropical plane cubic.}
\label{fig:1nodalcubic}
\end{figure}

\end{example}

The point is that each tropical plane degree-$d$ curve has a notion of
genus and a number of nodes, paralleling those
for algebraic plane curves. These notions let
us compute the Gromov-Witten invariant $N_{d, g}$ and the Severi
degree $N^{d, \delta}$ tropically, by Mikhalkin's Correspondence
Theorem (Theorem~\ref{thm:correspondencetheorem}) in
Section~\ref{sec:countingTropically}.

\subsection{Tropical Curves for Toric Surfaces}
\label{sec:toriccurves}

There is also a notion of tropical curves for any projective toric
surface (in this article all toric surfaces are assumed to be projective).

\begin{definition}
\label{def:tropicalcurvedegreeDelta}
A \emph{tropical curve of degree $\Delta$ } is a piecewise linear,
weighted graph $\Gamma$ in $\RR^2$ satisfying conditions
(\ref{itm:weights}), (\ref{itm:slopes}), and (\ref{itm:balanced}) of
Definition~\ref{def:tropicalcurve}, and additionally
\begin{enumerate}
\item[(\ref{itm:directions}')] the directions of the unbounded edges
  of $\Gamma$ are ``dual to $\partial \Lambda$``, i.e., 
the total
  weight of the unbounded edges of $\Gamma$ in direction $v$ equals
  the lattice length of the edge of $\Lambda$ with outer normal vector
  $v$.
\end{enumerate}
\end{definition}

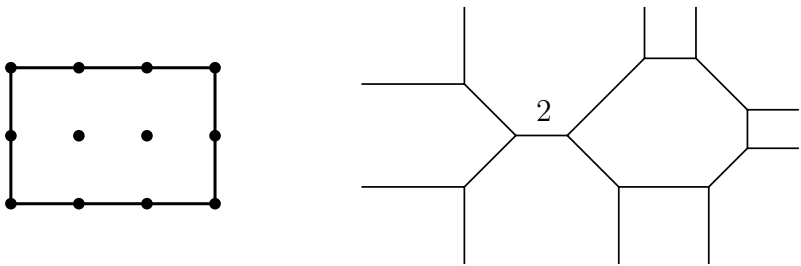
\begin{figure}[h]
\centering
\begin{minipage}{6in}
  \centering
  \raisebox{-0.5\height}{\includegraphics[height=2cm]{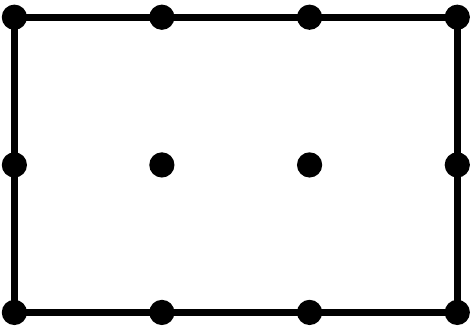}}
  \hspace*{.5in}
  \raisebox{-0.5\height}{
\begin{picture}(0,0)
\put(71,55){$2$}
\end{picture}
\includegraphics[height=3.5cm]{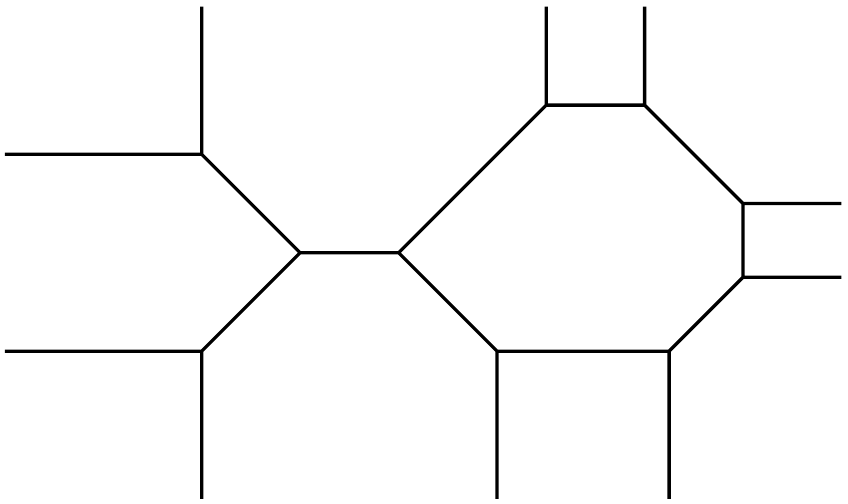}
}
\end{minipage}
\caption{Left: the Newton polygon $\Delta$ of curves in $\PP^1 \times \PP^1$ of
  bi-degree $(3,2)$. Right: a corresponding tropical curve of degree
  $\Delta$. All but one of its edges have weight $1$.}
\label{fig:P1timesP1tropicalcurve}
\end{figure}

Definition~\ref{def:tropicalcurvedegreeDelta} agrees with
Definition~\ref{def:tropicalcurve} for $\PP^2$ and degree $d$ curves:
in that case, $\Delta = \conv\{ (0,0), (d, 0), (0,d) \}$ defines the
toric surfaces $X(\Delta) = \PP^2$ with the ample line bundle
$\L(\Delta) = \O_{\PP^2}(d)$. The degree-$\Delta$ tropical curve has
$d$ rays in each of the directions $\binom{-1}{0}$, $\binom{0}{-1}$, and
$\binom{1}{1}$.

\begin{example}
A tropical curve of degree $\Delta$, with $\Delta = \conv\{ (0,0),
(3,0), (0,2), (3,2)\}$, is shown in
Figure~\ref{fig:P1timesP1tropicalcurve}. Such tropical curves count
algebraic curves in $\PP^1 \times \PP^1$ of bi-degree $(3,2)$, see Theorem~\ref{thm:toriccorrespondence}.
\end{example}

Just as in the $\PP^2$ case, we associate numerical invariants to a
tropical curve $\Gamma$ of degree $\Delta$. We say that $\Gamma$ is
\emph{irreducible}, if $\Gamma$ is not a union of two (non-empty)
tropical curves. The algebraic analog here is that the polynomial does not
factor. If
$\Gamma$ is irreducible, its \emph{genus} $g(\Gamma)$ is again defined
as the minimal first Betti number of any topological graph mapping
continuously and surjectively onto $\Gamma$. As before, if $\Gamma$
has irreducible components $\Gamma_1, \dots, \Gamma_r$, then its genus
is $\sum_{i=1}^r g(\Gamma_i) + 1 - r$.
The \emph{number of nodes} of a
tropical curve
$\Gamma$ is the number of interior lattice points $|\interior(\Delta) \cap
\ZZ^2|$ in $\Delta$ minus the genus $g(\Gamma)$ of $\Gamma$. This formula mimics
the algebraic analog of the relation between genus and the number of
nodes of a nodal algebraic curve in a toric surface.


\begin{remark}
\label{rmk:specialposition}
Despite allowing unbounded edges in
Definitions~\ref{def:tropicalcurve} and~\ref{def:tropicalcurvedegreeDelta} of weight
bigger than $1$, all weights of such edges are forced to be $1$ after
we fix sufficiently many point conditions (see
\cite[Lemma~4.20]{Mi03}).
Here, sufficiently many means that we require a finite curve count.
This will always be the case in
the sequel. Furthermore, by the same lemma, all vertices are
$3$-valent in such situations (where we view transversely crossing edges as
such).
\end{remark}

\section{Counting Algebraic Curves tropically}
\label{sec:countingTropically}

\subsection{Tropical Curve Enumeration for $\PP^2$}
\label{sec:tropEnumInP2}

Now we discuss how to enumerate algebraic curves
tropically, at least in some special cases. That this is indeed
possible can be seen as an instance of the ``shadow'' tropical
geometry containing enough information about the algebraic counterpart. We begin by considering the
complex projective plane $\PP^2$. Recall that the Gromov-Witten
invariant $N_{d, g}$ counts irreducible plane curves of degree $d$ and
genus $g$ passing through $3d + g - 1$ points in general position. The Severi
degree $N^{d, \delta}$ enumerates (not necessarily irreducible) plane
curves of degree $d$ with exactly $\delta$ nodes as singularities
through $\tfrac{(d+3)d}{2} - \delta$ points in general position.

It will turn out that, to enumerate algebraic curves with tropical curves,
we need to count them with a multiplicity. We define this multiplicity
for tropical curves with only $3$-valent vertices (and possibly some
transversely crossing edges). This will suffice to compute $N_{d, g}$ and $N^{d,
  \delta}$ by Remark~\ref{rmk:specialposition} as the tropical curves passing through sufficiently many
points in general position are automatically of this form.

Let $\Gamma$ be a tropical curve with $3$-valent vertex $v$. Let $e_1$
and $e_2$ be two of the $v$-adjacent edges of $\Gamma$. The
\emph{multiplicity} $\mult(v)$ of $v$ is
\[
\mult(v) = \wt(e_1) \cdot \wt(e_2) \cdot |\primitive(e_1, v) \wedge \primitive(e_2,
v)|,
\] 
thus $\mult(v)$ is $\wt(e_1) \cdot \wt(e_2)$ times the Euclidean area of the
parallelogram spanned by the two primitive vectors of $e_1$ and $e_2$. The balancing
condition implies the independence of which two edges we choose. The
\emph{multiplicity $\mult(\Gamma)$ of $\Gamma$} is the product over the
multiplicities of the $3$-valent vertices of $\Gamma$:
\begin{equation}
\label{eqn:multiplicity}
\mult(\Gamma) = \prod_{v \, 3\text{-valent}} \mult(v).
\end{equation}
In close analogy with the algebraic curve count, the \emph{tropical
  Gromov-Witten invariant} $N^{\trop}_{d, g}$ is the number of
irreducible tropical plane curves $\Gamma$ of degree $d$ and genus $g$ passing
through $3d+g-1$ points in $\RR^2$ in general position, counted with multiplicity
$\mult(\Gamma)$. Similarly, the
\emph{tropical Severi degree} $N_{\trop}^{d, \delta}$ is the number of
(possibly reducible) tropical plane curves $\Gamma$ of degree $d$ with $\delta$
nodes through $\tfrac{(d+3)d}{2} - \delta$ points in $\RR^2$ in general position,
counted with multiplicity $\mult(\Gamma)$. The following influential
theorem, suggested by Kontsevich and proved by Mikhalkin, says
that the tropical ``shadow'' suffices to compute the classical numbers
$N_{d, g}$ and $N^{d, \delta}$, and contributed a great deal to the
success of tropical geometry~\cite{AB10, FB, CJM10, CJM11, FM, GM052, GM053,MR08}.

\begin{theorem}[Mikhalkin's Correpondence Theorem for $\PP^2$
 {\cite[Theorem~1]{Mi03}}]
\hspace{2mm} 
\begin{enumerate}
\item We have $N_{d, g} = N^{\trop}_{d,g}$.
\item We have $N^{d, \delta} = N^{d, \delta}_{\trop}$.
\end{enumerate}
\label{thm:correspondencetheorem}
\end{theorem}

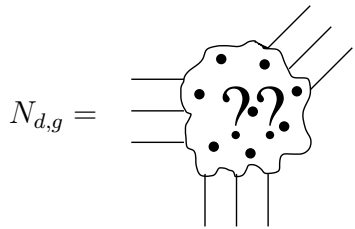
\begin{figure}[h]
\centering
\mbox{
\subfigure{
\text{
\begin{picture}(30,15)
\put(8,40){$N_{d, g} = $}
\end{picture}
}
}
\quad
\includegraphics[height=3cm]{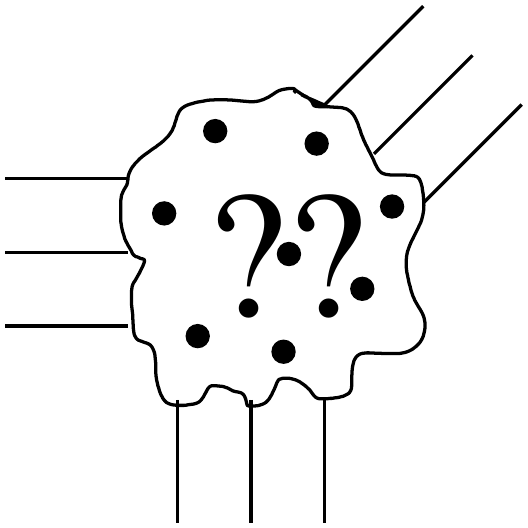}
}
\caption{Computing the Gromov-Witten invariant $N_{d, g}$ with
  piecewise linear geometry: $N_{d, g}$ equals the number of
  piecewise-linear tropical ``interpolation'' solutions, counted with
  multiplicity~(\ref{eqn:multiplicity}).}
\label{fig:cartoon}
\end{figure}

The Correspondence Theorem reduces the computation of $N_{d, g}$ and
$N^{d, \delta}$ to a piecewise-linear, combinatorial problem, which we
can think of as the ``cartoon'' in Figure~\ref{fig:cartoon}.

To compute $N_{3, 0}^{\trop}$, for example, we need to fill in the
picture in Figure~\ref{fig:cartoon} according to the rules of Definition~\ref{def:tropicalcurve},
with $3$ rays each going ``west'', ``south'', and ``north-east'',
respectively, together with $8$ point conditions. 

For larger $d$, this combinatorial problem is, however, still
subtle and difficult to carry out. A choice of a suitable point
configuration further simplifies the problem as the tropical curves
can be organized in a more manageable way. Specifically, if the points are
\emph{vertically stretched}, i.e., lie on a (classical) line of very
large, irrational slope (such a configuration is tropically generic), then the arising
tropical curves can be encoded by \emph{floor diagrams}, a family of decorated
graphs. For more details see~\cite{FB, BM2, FM}.

\begin{remark}
Mikhalkin shows a stronger statement than
Theorem~\ref{thm:correspondencetheorem}, namely an actual
``correspondence'' between algebraic curves through a particular
configuration and their tropicalizations.
On the complex $2$-torus $(\CC^*)^2$, we define, for any $t > 0$,
 \[
 \Log_t: 
 (\CC^*)^2 \to \RR^2, \quad (x,y) \mapsto (\log_t |x|, \log_t |y|).
 \]

Let $\P_\RR$ be a configuration of $3d+g-1$ tropically generic points
in $\RR^2$ (see \cite[Definition~4.7]{Mi03}). For $t > 0$, let
$\P_\CC^t$ be a configuration of $3d+g-1$ points in $\PP^2$ in general
position such that
$\Log_t(\P_\CC^t) = \P_\RR$. Let $\C_{\trop}$ be the set of irreducible tropical plane curves of
degree $d$ and genus $g$ through the points $\P_\RR$. Let $\C^t$ be
the set of irreducible complex plane curves of degree $d$ and genus
$g$ through the points $\C^t$. If we weight the tropical curves in $\C_{\trop}$ by their
multiplicity (\ref{eqn:multiplicity}), then, by
Theorem~\ref{thm:correspondencetheorem}, both sets $\C^t$ and 
$\C_{\trop}$ have cardinality $N_{d, g}$.

We can identify the curves in $\C^t$ and $\C_{\trop}$, for $t$ large
enough, as follows~\cite[Lemmas~8.3 and~8.4]{Mi03}: For each $\eps >0$, there is a $T > 0$ such that,
for $t \ge T$ and each tropical curve $\Gamma$ in $\C_{\trop}$, there are
precisely $\mult(\Gamma)$ complex curves $C$ in $\C^t$ with
\[
\Log_t(C) \subset \N_\eps(\Gamma) \subset \RR^2,
\]
where $\N_\eps(\Gamma)$ is an $\eps$-neighborhood of $\Gamma$.
Furthermore, each curve $C$ in $\C^t$ maps into $\N_\eps(\Gamma)$, for
some $\Gamma$ in $\C_{\trop}$. 

The content here is that
one can read off the cardinality of a fiber over $\Gamma$ of the tropicalization map
$\C^t \to \C_{\trop}$, for large $t$, from the tropical curve
$\Gamma$: it equals the tropical multiplicity $\mult(\Gamma)$.
\end{remark}

\begin{remark}
There is also a notion of ``parametrized'' tropical curves, which
are maps
$\pi: \Gamma^{\abs} \to \RR^2$ satisfying a balancing conditions,
from an ``abstract'' tropical curve $\Gamma^{\abs}$;
see~\cite[Section~2.2]{GM052} for the precise definition. In this language, the
tropical curves in Definitions~\ref{def:tropicalcurve}
and~\ref{def:tropicalcurvedegreeDelta} are the images
$\pi(\Gamma^{\abs})$. This notion is the tropical analog of stable maps
in Gromov-Witten theory and is, thus, more natural and flexible than embedded
curves in this setting.
In this paper, we chose to restrict to the simpler notion of
(embedded) tropical curves as in Definitions~\ref{def:tropicalcurve}
and~\ref{def:tropicalcurvedegreeDelta} as those are sufficient for our
purposes.
\end{remark}

\subsection{Tropical Curve Enumeration for Toric Surfaces}

A very similar approach works for arbitrary toric surfaces as
well. Recall that a lattice polygon $\Delta$ determines a projective
toric surface $X = X(\Delta)$, together with an ample line bundle
$\L = \L(\Delta)$. The number of irreducible genus-$g$ curves in the
complete linear system $|\L|$ passing through sufficiently many points
in general position is
denote $N_{\Delta, g}$; the number of (possibly reducible) curves in
$|\L|$ through $|\Lambda \cap \ZZ^2| - 1 - \delta$ points in general
position is
denoted by $N^{\Delta, \delta}$. The latter number $N^{\Delta,
  \delta}$ is known as the \emph{Severi degree of the 
  surface $X(\Delta)$}.

As before, we define a tropical analog of the numbers $N_{\Delta, g}$
and $N^{\Delta, \delta}$. Let $N_{\Delta, g}^{\trop}$ be the number of
irreducible tropical degree-$\Delta$ curves $\Gamma$ with genus $g$
through sufficiently many points in $\RR^2$ in general position,
counted with multiplicity $\mult(\Gamma)$ (see
(\ref{eqn:multiplicity})). Let $N^{\Delta, \delta}_{\trop}$ be the
number of (possibly reducible) tropical degree-$\Delta$ curves
$\Gamma$ with $\delta$ nodes through $|\Delta \cap \ZZ^2| - 1 -
\delta$ points in $\RR^2$ in general position, counted with multiplicity $\mult(\Gamma)$.

\begin{theorem}[Correspondence Theorem for Toric Surfaces
  {\cite[Theorem~1]{Mi03}}]
\label{thm:toriccorrespondence}
\hspace{30mm}
\begin{enumerate}
\item We have $N_{\Delta, g} = N_{\Delta, g}^{\trop}$.
\item We have $N^{\Delta, \delta} = N^{\Delta, \delta}_{\trop}$.
\end{enumerate}
\end{theorem}

\section{Applications}
\label{sec:applications}
While it's certainly nice to have
combinatorial descriptions of classical curve enumeration problems (as
in Section~\ref{sec:countingTropically}),
the power of tropical techniques comes with their ability to
prove deep and new theorems in enumerative algebraic geometry. In this
section, we collect a few of these applications. For more,
see for example ~\cite{BGM10, GM052, GM053}.

\subsection{Node Polynomials for Plane Curves}


Steiner~\cite{St48},
in 1848, computed the degree $N^{d, 1} = 3(d-1)^2$ of the discriminant
of $\PP^2$. A few decades later, in 1863 resp.\ 1867, Cayley resp.\
Roberts, gave polynomial expressions for $N^{d, 2}$
resp.\ $N^{d, 3}$ (in the latter case for $d \ge 3$).

Much later, in 1994, Di Francesco and Itzykson~\cite{DI}, conjectured the numbers
$N^{d, \delta}$ to be polynomial in $d$, for fixed $\delta$ and $d$ large enough.
For $\delta = 4, 5$, and $6$, this was affirmed by
Vainsencher \cite{Va} in 1995 using deformation theory. In 2001, Kleiman and Piene
\cite{KP} settled the cases $\delta = 7, 8$ utilizing similar
techniques.
Fomin and Mikhalkin~\cite{FM}, in 2009,  proved Di Francesco and
Itzykson's conjecture, using tropical geometry techniques.

\begin{theorem} [{\cite[Theorem~5.1]{FM}}]
\label{thm:planarnodepolys}
For fixed $\delta \ge 1$, there is combinatorially defined polynomial
$N_\delta(d)$ in $d$, such that
\[
N_{\delta}(d) = N^{d, \delta},
\]
provided that $d \ge 2 \delta$.
\end{theorem}

Here ``combinatorially defined'' means that Fomin and Mikhalkin's
description of the polynomials $N_\delta(d)$ gives rise to
a combinatorial algorithm computing $N_\delta(d)$. Their method was
improved and implemented by the author~\cite{FB}, who computed
$N_\delta(d)$ for $\delta \le 14$.
 Following Kleiman and Piene~\cite{KP}, the $N_\delta(d)$ are 
called \emph{node polynomials}.

Fomin and Mikhalkin's proof is mostly combinatorial and uses a
description of $N^{d, \delta}$ in terms of \emph{floor diagrams}. These
purely combinatorial
gadgets, introduced
by Brugall\'e and Mikhalkin~\cite{BM1, BM2}, are a
family of enriched graphs, arising from tropical plane curves by
topological contractions (see~\cite[Section~3]{FM}
or \cite[Section~4]{BM2}).

There are now also alternate proofs of
Theorem~\ref{thm:planarnodepolys}.  Tzeng, in her celebrated
 work~\cite{Tz10}, proved the G\"ottsche conjecture~\cite{Go} using
 algebraic cobordism. (For the precise statement of the conjecture,
 see Section~\ref{sec:toricpolys} below.) Tzeng thus
 established universal polynomiality (in certain Chern numbers) of
 the Severi degree of any smooth projective surface.
A second proof of the G\"ottsche conjecture was given soon after by Kool,
Shende, and Thomas~\cite{KST11} using BPS calculus~\cite{PT10}. Although the
algebro-geometric techniques give rise to an algorithm to compute node
polynomials, the tropical approach seems to be the most efficient, at
least in the case of $\PP^2$ and $\PP^1 \times \PP^1$. By the
universality of the node polynomials, the two 
cases suffice to determine the Severi degree of \emph{any} smooth projective
surface (and the two cases can indeed be computed by tropical geometry).

\subsection{Node Polynomials for Toric Surfaces}
\label{sec:toricpolys}

Tropical geometry techniques can also be used to compute node
polynomials for a large family of toric surfaces that are,  in
general, non-smooth. Recall that a lattice polygon $\Delta$ determines
a toric surface $X(\Delta)$, together with a line bundle
$\L(\Delta)$. Given $\Delta$, we are interested in the number of
$\delta$-nodal curves in the complete linear system $|\L|$ through
$|\Delta \cap \ZZ^2| - 1 - \delta$ points in general position. As mentioned before, this number is
the Severi degree $N^{\Delta, \delta}$. This section is based
on~\cite{AB10}.

In the following, we restrict the presentation to the special case of
Hirzebruch surfaces in favor of simpler notation. Such surfaces are
smooth, but exhibit (see Theorem~\ref{thm:toricnodepolynomials})
already the main features of the more general case: the Severi degrees
are polynomial in the ``multi-degree'' $\L(\Delta)$ and parameters
of the surface $X(\Delta)$. 

\begin{figure}[h]
\centering
\mbox{
\subfigure{
\includegraphics[height=2.0cm]{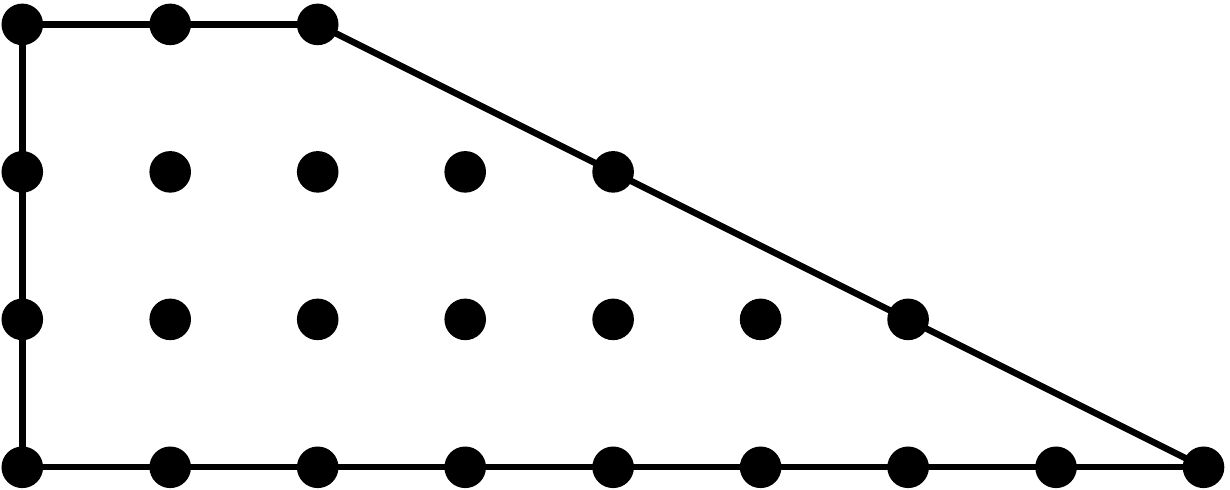}
}
}
\caption{
The polygon $\Delta$ of the Hirzebruch surfaces $F_2$ and a divisor of
  bi-degree $(2, 3)$.
}
\label{fig:hirzebruchpolygon}
\end{figure}

Let $N^{(a,b), \delta}_m$ be the number
of $\delta$-nodal curves in the linear system determined by a divisor 
of \emph{bi-degree} $(a, b)$ on the Hirzebruch surface $F_m$, i.e.,
with Newton polygon $\conv((0,0), (0,b), (a,b), (a+bm,0))$, up to translation. For $m = 0$,
this means enumerating the $\delta$-nodal curves in $\PP^1 \times
\PP^1$ of bi-degree $(a, b)$ through $(a+1)(b+1) - 1 - \delta$
points in general position. The polygon determining $F_m$ and the corresponding line
bundle is shown in Figure~\ref{fig:hirzebruchpolygon}.

\begin{theorem}
\label{thm:toricnodepolynomials}
For every $\delta \ge 1$, there is a combinatorially defined polynomial $p_\delta(m, a,
b)$ such that $N^{(a, b), \delta}_m = p_\delta(m, a, b)$ provided $a +
m \ge 2 \delta $ and $b \ge 2 \delta$.
\end{theorem}

G\"ottsche~\cite[Conjecture~2.1]{Go} famously conjectured the existence of
universal polynomials $T_\delta(x, y, z, t)$ that compute the Severi degree for any smooth projective
surface $S$ and any sufficiently ample line bundle $\L$ on $S$. According to the
conjecture, the number of
$\delta$-nodal curves in the linear system $|\L|$ through an
appropriate number of points is given by evaluating $T_\delta$ at
the four topological numbers 
$\L^2, \L K_S, K_S^2$ and $c_2(S)$. 
Here, $K_S$ denotes the
canonical bundle, $c_1$ and $c_2$ represent Chern classes, and $LM$
denotes the degree of $c_1(L) \cdot c_1(M)$ for line bundles $L$ and
$M$.
In the setting of Theorem~\ref{thm:toricnodepolynomials}, the
four topological numbers are polynomial in $m$, $a$, and $b$. The theorem thus
also follows from the G\"ottsche conjecture.

One can prove a similar result as in
Theorem~\ref{thm:toricnodepolynomials} for ``h-transverse'' polygons
$\Delta$ \cite[Theorem~1.3]{AB10}. Such polygons are allowed to have only edges of slope $1/n$, for $n
\in \ZZ \cup \{ \infty \}$. The resulting toric surfaces are not
smooth in general and are, thus, outside the realm of the G\"ottsche
conjecture. 
Still, the Severi degree is polynomial in
parameters of $\Delta$, and we can use tropical geometry to prove it.

\subsection{Double Hurwitz Numbers}

Fix two partitions $\lambda = (\lambda_1, \lambda_2, \dots, \lambda_m>0)$ and $\mu
= (\mu_1, \mu_2, \dots, \mu_n>0)$ of a positive integer $d$. The \emph{double Hurwitz
  number} $H_g(\lambda; \mu)$ counts degree-$d$ maps $\pi: C \to
\PP^1$, where $C$ is a connected, genus $g$ curve and $\pi$ has
\emph{ramification profiles} $\lambda$ resp.\ $\mu$ over $0$ resp.\
$\infty$, and simple ramification over $r = 2g - 2 + m + n$ fixed
other points. Each cover is counted with weight $1/|\Aut(\pi)|$.
The main reference for this section is Cavalieri, Johnson, and
Markwig's paper~\cite{CJM11}.

We can think of double Hurwitz numbers as the $\PP^1$-analog of the
Gromov-Witten invariants $N_{d, g}$ of $\PP^2$: instead of counting
maps to $\PP^2$ (each degree-$d$ plane curve of genus $g$ is the image of a
degree-$d$ map from
an abstract genus-$g$ curve to $\PP^2$), we now count such maps to
$\PP^1$.

Alternatively, $H_g(\lambda; \mu)$ counts tuples of permutations
$\sigma_0, \sigma_1, \dots, \sigma_r, \sigma_\infty \in S_d$ with
\smallskip
\begin{itemize}
\item $\sigma_0$ and $\sigma_\infty$ have cycle type $\lambda$ and
  $\mu$, respectively,
\item $\sigma_1, \dots, \sigma_r$ are transpositions,
\item $\sigma_0 \sigma_1 \cdots \sigma_r \sigma_\infty = \id \in S_d$,
\item the subgroup generated by $\sigma_0, \dots, \sigma_\infty$ acts
  transitively on $\{1, \dots, d\}$.
\end{itemize}
\smallskip
We weight the count by $d! \cdot |\Aut(\sigma_0)| \cdot
|\Aut(\sigma_\infty)|$. Here, the number $r$ of points with simple
ramification is determined by the Riemann-Hurwitz formula: $r = 2g - 2
+ |\lambda| + |\mu|$.

Goulden, Jackson, and Vakil~\cite{GJV05} showed that, for $\lambda$
and $\mu$ of fixed length, double Hurwitz
numbers $H_g(\lambda; \mu)$ are piecewise polynomial in the entries of
$\mu$ and $\nu$. This means that there is a hyperplane arrangement in
$\RR^{l(\lambda) + l(\mu)}$ such that  $H_g(\lambda; \mu)$ is
polynomial on each connected component of the complement of the
arrangement. Here, $l(\lambda)$ and $l(\mu)$ are the number of parts
of $\lambda$ and $\mu$.
Using a tropical analog of Hurwitz numbers, Cavalieri, Johnson, and
Markwig~\cite{CJM10} confirmed Goulden, Jackson, and Vakil's result.

\begin{theorem}[{\cite[Theorem 2.1]{GJV05}}, {\cite[Theorem~1.1]{CJM11}}]
For fixed $l(\lambda)$ and $l(\mu)$, the function $H_g(\lambda; \mu)$: $(\ZZ-\{ 0 \})^{\l(\lambda) +
  l(\mu)} \to \QQ$ is piecewise polynomial.
\end{theorem}

Shadrin, Shapiro, and Vainshtein~\cite{SSV08} computed, in genus $0$, the
chamber structure of $H_0(\lambda; \mu)$ (i.e., the domains of polynomiality) as well as an explicit
``wall crossing formula.'' The latter describes how $H_0(\lambda; \mu)$
changes when one moves from a chamber to an adjacent chamber.

Cavalieri, Johnson, and Markwig~\cite{CJM11} generalized these
results to double Hurwitz numbers for all genera, using tropical
geometry techniques, and thus giving a common approach to
both~\cite{GJV05} and~\cite{SSV08}.

 \begin{theorem}[{\cite[Theorems~1.3 and~1.5]{CJM11}}]
 \textcolor{white}{ }
 \begin{enumerate}
 \item The chambers of polynomiality of $H_g(\lambda;\mu)$ are the
   complements of hyperplanes given by explicit formulas.
 \item The wall crossing formulas are computed explicitly in terms
   of Hurwitz numbers with fewer simple ramification points.
\end{enumerate}
 \end{theorem}

For the explicit formulas, see~\cite{CJM11}.



\subsection{Real Enumerative Geometry via Tropical Geometry}
\label{sec:Welschinger}

Tropical geometry was also successfully applied to problems in real
algebraic geometry. There, one studies, for example, those complex solutions to polynomial
equation that are invariant (as a set) under complex conjugation. In
general, the real analogs of complex enumerative problems are even more
difficult. Nevertheless, some of them can be addressed nicely by
tropical means. For a general reference on real algebraic geometry,
see~\cite{So03}, for some further tropical application see for example \cite{IKS03, IKS09,Mi03}.

In this section, we focus on enumeration of \emph{real plane curves},
i.e., complex algebraic curves in $\PP^2$ invariant under complex
conjugation. A natural question is about a real analog of
Gromov-Witten invariants of $\PP^2$: real plane curves of fixed degree
and genus passing through a real point configuration in general
position. We quickly run into difficulty
however: if we try to naively count such curves, we find that their number does
depend on the point configuration! For example, the number
of real rational cubics through $8$ real points in general position can be
$8$, $10$, or $12$.

Welschinger~\cite{We03} resolved this problem by proposing to count
real plane curves with a sign. To a real plane curve $C$, Welschinger associated a
multiplicity $(-1)^s$, where $s$ is the number of real double points
of $C$ (a real double point is locally given by $\{ x^2 + y^2 = 0\}$,
with $x$ and $y$ some local
coordinates). 
For $d \ge 1$, let the \emph{Welschinger invariant} $W_d$ be
the number of rational degree-$d$ real plane curves, counted with
Welschinger's multiplicity, passing through $3d-1$ points in general
position.  Welschinger showed that $W_d$ is indeed independent of the point
configuration \cite[Theorem~2.1]{We03}.
 In particular, there are always $W_3 = 8$ real rational
plane cubics through $8$ points in general position, if one counts them
with Welschinger's multiplicity. 

From the definition, it is not at all obvious whether $W_d$ is positive
or negative or zero. With tropical geometry, Itenberg,
Kharlamov, and Shustin resolved this. To the author's
knowledge, no non-tropical proof of the following has been discovered yet.

\begin{theorem}[{\cite[Theorem~1.1]{IKS03}}]
For all $d \ge 1$, the Welschinger invariant $W_d$ is positive.
\end{theorem}

The proof is based on a real analog of Mikhalkin's
Correspondence Theorem: Mikhalkin associates in~\cite{Mi03} to a
tropical curve $\Gamma$ not only a (complex) multiplicity
$\mult(\Lambda)$ (as we do in Section~\ref{sec:tropEnumInP2}), but also a
\emph{real multiplicity} $\mult_\RR(\Lambda)$. For the precise
definition, see \cite[Definition~7.19]{Mi03}. Similarly to the complex case,
we define the \emph{tropical Welschinger invariant} $W_d^{\trop}$ as
the number of tropical degree-$d$ genus-$0$ plane curves $\Gamma$ passing
though $3d-1$ points in general position, but now counted with multiplicity $\mult_\RR(\Gamma)$.

\begin{theorem}
[{Mikhalkin's Real Correspondence Theorem \cite[Theorem 6]{Mi03}}]
\textcolor{white}{ } \newline
For any $d \ge 1$, we have
\[
W_d = W_d^{\trop}.
\]
\end{theorem}

\bibliographystyle{amsplain}
\bibliography{References_Florian}

\end{document}